\documentstyle[12pt,pb-diagram]{article}   
\newcommand{\be}{\begin{equation}}
\newcommand{\ee}{\end{equation}}
\newcommand{\bea}{\begin{eqnarray}}
\newcommand{\eea}{\end{eqnarray}}
\newcommand{\bean}{\begin{eqnarray*}}
\newcommand{\eean}{\end{eqnarray*}}
\newcommand{\N}{I\!\!N}
\newcommand{\R}{I\!\!R}
\newcommand{\Z}{Z\!\!\!Z}
\newcommand{\PP}{I\!\!P}

\newcommand{\C}{\,I\!\!\!\!C}
\newcommand{\Q}{I\!\!\!\!Q}

\newcommand{\CE}{{\cal E}}

\newcommand{\CB}{{\cal B}}
\newcommand{\CC}{{\cal C}}

\newcommand{\CD}{{\cal D}}

\newcommand{\Spec}{\mbox{Spec}\;}

\newcommand{\supp}{\mbox{supp}\,}

\newcommand{\veee}{\scriptscriptstyle\vee}

\newcounter{Abschnitt}[section]
\newcommand{\neu}[1]{\protect\refstepcounter{Abschnitt}\protect\label{t#1}\vspace{1ex}
                     {\bf (\arabic{section}.\protect\arabic{Abschnitt})}
                     $\qquad$}
\newcommand{\zitat}[2]{(\protect\ref{s#1}.\protect\ref{t#1#2})}
\newcommand{\surj}{\longrightarrow\hspace{-1.5em}\longrightarrow}

\newcommand{\kgeq}{\geq}
\newcommand{\ku}{\underline}
\newcommand{\ks}{\scriptstyle}
\newcommand{\kd}{\displaystyle}

\newcommand{\kR}{{\cal R}}

\newcounter{secnum}

\newcommand{\sect}[1]
 {\protect\section{#1}
  \protect\setcounter{secnum}{\value{section}}
  \protect\setcounter{equation}{0}
 \protect\renewcommand{\theequation}{\mbox{\arabic{secnum}.\arabic{equation}}}}

\begin{document}
\title{Infinitesimal Deformations and Obstructions for Toric Singularities}

\author{Klaus Altmann\footnotemark[1]\\
        \small Dept.~of Mathematics, M.I.T., Cambridge, MA 02139, U.S.A.
        \vspace{-0.7ex}\\ \small E-mail: altmann@math.mit.edu}
\footnotetext[1]{Die Arbeit wurde mit einem Stipendium des DAAD unterst\"utzt.}
\date{}
\maketitle

\begin{abstract}
The obstruction space $T^2$ and the cup product
$T^1\times T^1\to T^2$ are computed for toric singularities.
\end{abstract}

\tableofcontents

%
%
\sect{Introduction}\label{s1}


\neu{11}
For an affine scheme $\,Y= \Spec A$, there are two important $A$-modules,
$T^1_Y$ and $T^2_Y$, carrying information about its deformation theory:
$T^1_Y$ describes the infinitesimal deformations, and $T^2_Y$ contains the
obstructions for extending deformations of $Y$ to larger base spaces.\\
\par
In case $Y$ admits a versal deformation, $T^1_Y$ is the tangent space of the
versal base space $S$. Moreover, if $J$ denotes the ideal defining $S$ as a
closed
subscheme of the affine space $T^1_Y$, the module
$\left( ^{\displaystyle J}\! / \! _{\displaystyle m_{T^1} \,J} \right) ^\ast$
can be canonically embedded into $T^2_Y$, i.e. $(T_Y^2)^\ast$-elements induce
the
equations defining $S$ in $T^1_Y$.\\
\par
The vector spaces $T^i_Y$ come with a cup product
$T_Y^1 \times T^1_Y \rightarrow T^2_Y$.
The associated quadratic form $T^1_Y \rightarrow T^2_Y$ describes the
quadratic part of the elements of $J$, i.e. it can be used to get a better
approximation of the versal base space $S$ as regarding its tangent space
only.\\
\par


\neu{12}
In \cite{T1} we have determined the vector space $T^1_Y$ for affine toric
varieties.
The present paper can be regarded as its continuation - we will compute $T^2_Y$
and
the cup product. \\
These modules $T^i_Y$ are canonically graded (induced from the character group
of
the torus). We will describe their homogeneous pieces as cohomology groups of
certain complexes, that are directly induced from the combinatorial structure
of the
rational, polyhedral cone defining our variety $Y$.  The results
can be found in \S \ref{s3}.\\
\par
Switching to another, quasiisomorphic complex provides a second formula for the
vector spaces $T^i_Y$ (cf. \S \ref{s6}). We will use this particular version
for
describing these spaces and the cup product in the special case of
three-dimensional toric Gorenstein singularities (cf. \S \ref{s7}).\\
\par

%
%
\sect{$T^1$, $T^2$, and the cup product (in general)}\label{s2}

In this section we will give a brief reminder to the well known
definitions of these objects. Moreover, we will use this opportunity to fix
some
notations.\\
\par


\neu{21}
Let $Y \subseteq \C^{w+1}$ be given by equations $f_1,\dots,f_m$, i.e.
its ring of regular functions equals
\[
A=\;^{\displaystyle P}\!\! / \! _{\displaystyle I} \quad
\mbox{ with }
\begin{array}[t]{l}
P = \C[z_0,\dots, z_w]\\
I = (f_1,\dots,f_m)\, .
\end{array}
\]
Then, using $d:^{\displaystyle I}\! / \! _{\displaystyle I^2}
\rightarrow A^{w+1}\;$ ($d(f_i):= (\frac{\partial f_i}{\partial z_0},\dots
\frac{\partial f_i}{\partial z_w})$),
the vector space $T^1_Y$ equals
\[
T^1_Y = \;^{\displaystyle \mbox{Hom}_A(^{\displaystyle I}\!\! / \!
_{\displaystyle I^2}, A)} \! \left/ \!
_{\displaystyle \mbox{Hom}_A(A^{w+1},A)} \right.\; .
\vspace{1ex}
\]
\par


\neu{22}
Let $\kR\subseteq P^m$ denote the $P$-module of relations between the equations
$f_1,\dots,f_m$. It contains the so-called Koszul relations
$\kR_0:= \langle f_i\,e^j - f_j \,e^i \rangle$ as a submodule.\\
Now, $T^2_Y$ can be obtained as
\[
T^2_Y = \;^{\displaystyle \mbox{Hom}_P(^{\displaystyle \kR}\! / \!
_{\displaystyle \kR_0}, A)} \! \left/ \!
_{\displaystyle \mbox{Hom}_P(P^m,A)} \right.\; .
\vspace{1ex}
\]
\par


\neu{23}
Finally, the cup product $T^1\times T^1 \rightarrow T^2$ can be defined in the
following way:
\begin{itemize}
\item[(i)]
Starting with an $\varphi\in \mbox{Hom}_A(^{\displaystyle I}\! / \!
_{\displaystyle I^2}, A)$, we lift the images of the $f_i$ obtaining
elements $\tilde{\varphi}(f_i)\in P$.
\item[(ii)]
Given a relation $r\in \kR$, the linear combination
$\sum_ir_i\,\tilde{\varphi}(f_i)$ vanishes in $A$, i.e. it is contained in the
ideal $I\subseteq P$. Denote by $\lambda(\varphi)\in P^m$ any set of
coefficients such that
\[
\sum_i r_i \, \tilde{\varphi}(f_i) + \sum_i \lambda_i(\varphi)\, f_i =0\quad
\mbox{ in } P.
\]
(Of course, $\lambda$ depends on $r$ also.)
\item[(iii)]
If $\varphi, \psi \in \mbox{Hom}_A(^{\displaystyle I}\! / \!
_{\displaystyle I^2}, A)$ represent two elements of $T^1_Y$, then we define for
each relation $r\in \kR$
\[
(\varphi \cup \psi)(r) := \sum_i \lambda_i (\varphi)\, \psi(f_i) +
\sum_i \varphi(f_i)\, \lambda_i(\psi)\; \in A\, .
\vspace{1ex}
\]
\end{itemize}
{\bf Remark:}
The definition of the cup product does not depend on the choices we made:
\begin{itemize}
\item[(a)]
Choosing a $\lambda'(\varphi)$ instead of $\lambda(\varphi)$ yields
$\lambda'(\varphi) - \lambda(\varphi) \in \kR$, i.e. in $A$ we obtain the same
result.
\item[(b)]
Let $\tilde{\varphi}'(f_i)$ be different liftings to $P$. Then, the difference
$\tilde{\varphi}'(f_i) - \tilde{\varphi}(f_i)$ is contained in $I$, i.e. it can
be written as some linear combination
\[
\tilde{\varphi}'(f_i) - \tilde{\varphi}(f_i) = \sum_j t_{ij}\, f_j\, .
\]
Hence,
\[
\sum_i r_i \,\tilde{\varphi}'(f_i) = \sum_i r_i \,\tilde{\varphi}(f_i) +
\sum_{i,j} t_{ij}\, r_i\, f_j\,,
\]
and we can define $\lambda'_j(\varphi):= \lambda_j(\varphi) -
\sum_it_{ij}\,r_i$
(corresponding to $\tilde{\varphi}'$ instead of $\tilde{\varphi}$). Then, we
obtain
for the cup product
\[
(\varphi\cup\psi)'(r) - (\varphi\cup\psi)(r) = -\sum_ir_i\cdot
\left( \sum_j t_{ij}\, \psi(f_j)\right)\, ,
\]
but this expression comes from some map $P^m\rightarrow A$.
\vspace{3ex}
\end{itemize}

%
%
\sect{$T^1$, $T^2$, and the cup product (for toric varieties)}\label{s3}


\neu{31}
We start with fixing the usual notations when dealing with affine toric
varieties (cf. \cite{Ke} or
\cite{Oda}):
\begin{itemize}
\item
Let $M,N$ be mutually dual free Abelian groups, we denote by $M_{\R}, N_{\R}$
the associated real
vector spaces obtained by base change with $\R$.
\item
Let $\sigma=\langle a^1,\dots,a^N\rangle \subseteq N_{\R}$ be a rational,
polyhedral
cone with apex - given by its fundamental generators. \\
$\sigma^{\veee}:= \{ r\in M_{\R}\,|\; \langle \sigma,\,r\rangle \geq 0\}
\subseteq M_{\R}$
is called the dual cone. It induces a partial order on the lattice $M$ via
$[\,a\kgeq b$ iff
$a-b \in \sigma^{\veee}\,]$.
\item
$A:= \C[\sigma^{\veee}\cap M]$ denotes the semigroup algebra. It is the ring of
regular
functions of the toric variety $Y= \Spec A$ associated to $\sigma$.
\item
Denote by $E\subset \sigma^{\veee}\cap M$ the minimal set of generators of this
semigroup
("the Hilbert basis"). $E$ equals the set of all primitive (i.e.
non-splittable) elements
of $\sigma^{\veee}\cap M$.
In particular, there is a surjection of semigroups $\pi:\N^E \surj
\sigma^{\veee}\cap M$, and
this fact translates into a closed embedding $Y\hookrightarrow \C^E$.\\
To make the notations
coherent with \S \ref{s2}, assume that $E=\{r^0,\dots,r^w\}$ consists of $w+1$
elements.
\vspace{2ex}
\end{itemize}


\neu{32}
To a fixed degree $R\in M$ we associate ``thick facets'' $K_i^R$ of the dual
cone
\[
K_i^R := \{r\in \sigma^{\veee}\cap M \, | \; \langle a^i, r \rangle <
\langle a^i, R \rangle \}\quad (i=1,\dots,N) .
\vspace{2ex}
\]
\par

{\bf Lemma:}{\em
\begin{itemize}
\item[(1)]
$\cup_i K_i^R = (\sigma^{\veee}\cap M) \setminus (R+ \sigma^{\veee})$.
\item[(2)]
For each $r,s\in K_i^R$ there exists an $\ell\in K_i^R$ such that $\ell\geq
r,s$.
Moreover, if $Y$ is smooth in codimension 2, the intersections $K^R_i\cap
K^R_j$
(for 2-faces $\langle a^i,a^j\rangle <\sigma$) have the same property.
\vspace{1ex}
\end{itemize}
}
\par

{\bf Proof:}
Part (i) is trivial; for (ii) cf. (3.7) of \cite{T1}.
\hfill$\Box$\\
\par
Intersecting these sets with $E\subseteq \sigma^{\veee}\cap M$, we obtain the
basic objects for describing the modules $T^i_Y$:
\bean
E_i^R &:=& K_i^R \cap E = \{r\in E\, | \; \langle a^i,r \rangle <
\langle a^i, R \rangle \}\, ,\\
E_0^R &:=& \bigcup_{i=1}^N E_i^R\; ,\mbox{ and}\\
E^R_{\tau} &:=& \bigcap_{a^i\in \tau} E^R_i \; \mbox{ for faces }
\tau < \sigma\,.
\eean
We obtain a complex
$L(E^R)_{\bullet}$ of free Abelian groups via
\[
L(E^R)_{-k} := \bigoplus_{\begin{array}{c}
\tau<\sigma\\ \mbox{dim}\, \tau=k \end{array}} \!\!L(E^R_{\tau})
\]
with the usual differentials.
($L(\dots)$ denotes the free Abelian group of integral, linear dependencies.)
\\
The most interesting part ($k\leq 2$) can be written explicitely as
\[
L(E^R)_{\bullet}:\quad \cdots
\rightarrow
\oplus_{\langle a^i,a^j\rangle<\sigma} L(E^R_i\cap E^R_j)
\longrightarrow
\oplus_i L(E_i^R) \longrightarrow L(E_0^R) \rightarrow 0\,.
\vspace{1ex}
\]
\par


\neu{33}
{\bf Theorem:}
{\em
\begin{itemize}
\item[(1)]
$T^1_Y(-R) = H^0 \left( L(E^R)_{\bullet}^\ast \otimes_{\Z}\C\right)$
\item[(2)]
$T^2_Y(-R) \supseteq H^1 \left( L(E^R)_{\bullet}^\ast \otimes_{\Z}\C\right)$
\item[(3)]
Moreover, if $Y$ is smooth in codimension 2
(i.e.\ if the 2-faces $\langle a^i, a^j \rangle < \sigma$ are spanned
by a part of a $\Z$-basis of the lattice $N$), then
\[
T^2_Y(-R) = H^1 \left( L(E^R)_{\bullet}^\ast \otimes_{\Z}\C\right)\, .
\]
\item[(4)]
Module structure: If $x^s\in A$ (i.e. $s\in \sigma^{\veee}\cap M$), then
$E_i^{R-s}\subseteq R_i^R$, hence $L(E^R)_{\bullet}^\ast \subseteq
L(E^{R-s})^\ast_{\bullet}$. The induced map in cohomology corresponds to the
multiplication with $x^s$ in the $A$-modules $T^1_Y$ and $T^2_Y$.
\vspace{2ex}
\end{itemize}
}
\par

The first part was shown in \cite{T1}; the formula for $T^2$ will be proved
in \S \ref{s4}. Then, the claim concerning the module structure will become
clear
automatically.\\
\par

{\bf Remark:}
The assumption made in (3) can not be dropped: \\
Taking for $Y$ a 2-dimensional cyclic quotient
singularity given by some 2-dimensional cone $\sigma$, there are only two
different sets
$E_1^R$ and $E_2^R$ (for each $R\in M$). In particular,
$H^1 \left( L(E^R)_{\bullet}^\ast \otimes_{\Z}\C\right)=0$.\\
\par


\neu{34}
We want to describe the isomorphisms connecting the general $T^i$-formulas of
\zitat{2}{1} and \zitat{2}{2} with the toric ones given in \zitat{3}{3}.\\
\par
$Y\hookrightarrow\C^{w+1}$ is given by the equations
\[
f_{ab}:= \ku{z}^a-\ku{z}^b\quad (a,b\in \N^{w+1} \mbox{ with } \pi(a)=\pi(b)
\mbox{ in } \sigma^{\veee}
\cap M),
\]
and it is easier to deal with this infinite set of equations
(which generates the ideal $I$ as a $\C$-vector
space) instead of selecting a finite number of them in some non-canonical way.
In particular, for
$m$ of \zitat{2}{1} and \zitat{2}{2} we take
\[
m:= \{ (a,b)\in \N^{w+1}\times\N^{w+1}\,|\;\pi(a)=\pi(b)\}\,.
\]
The general $T^i$-formulas mentioned in \zitat{2}{1} and \zitat{2}{2} remain
true.\\
\par

{\bf Theorem:}
{\em
For a fixed element $R\in M$ let $\varphi: L(E)_{\C}\rightarrow \C$ induce some
element of
\[
\left(\left. ^{\displaystyle L(E_0^R)}\!\!\right/
\!_{\displaystyle \sum_i L(E_i^R)} \right)^\ast \otimes_{\Z} \C
\cong T^1_Y(-R)\quad \mbox{(cf. Theorem \zitat{3}{3}(1)).}
\]
Then, the $A$-linear map
\bean
^{\displaystyle I}\!\!/\!_{\displaystyle I^2} &\longrightarrow& A\\
\ku{z}^a-\ku{z}^b & \mapsto & \varphi(a-b)\cdot x^{\pi(a)-R}
\eean
provides the same element via the formula \zitat{2}{1}.
}\\
\par
Again, this Theorem follows from the paper \cite{T1} - accompanied by the
commutative diagram
of \zitat{4}{3} in the present paper. (Cf. Remark \zitat{4}{4}.)\\
\par

{\bf Remark:}
A simple, but nevertheless important check shows that the map
$(\ku{z}^a-\ku{z}^b) \mapsto
\varphi(a-b)\cdot x^{\pi(a)-R}$ goes into $A$, indeed:\\
Assume $\pi(a)-R \notin \sigma^{\veee}$. Then, there exists an index $i$ such
that
$\langle a^i, \pi(a)-R \rangle <0$.
Denoting by "supp $q$" (for a $q\in \R^E$) the set of those $r\in E$ providing
a non-vanishing entry
$q_r$, we obtain
\[
\mbox{supp}\,(a-b) \subseteq \mbox{supp}\,a \cup \mbox{supp}\, b \subseteq
E^R_i\, ,
\]
i.e. $\varphi(a-b)=0$.\\
\par


\neu{35}
The $P$-module $\kR\subseteq P^m$ is generated by relations of two different
types:
\bean
r(a,b;c) &:=& e^{a+c,\,b+c}- \ku{z}^c\, e^{a,b}\quad
(a,b,c\in \N^{w+1};\, \pi(a)=\pi(b))\quad \mbox{ and}\\
s(a,b,c) &:=& e^{b,c} - e^{a,c} + e^{a,b}\quad
(a,b,c\in \N^{w+1};\, \pi(a)=\pi(b)=\pi(c))\,.\\
&&\qquad\qquad(e^{\bullet,\bullet} \mbox{ denote the standard basis vectors of
} P^m.)
\vspace{1ex}
\eean
\par


{\bf Theorem:}
{\em
For a fixed element $R\in M$ let $\psi_i: L(E_i^R)_{\C}\rightarrow \C$ induce
some
element of
\[
\left( \frac{\displaystyle
\mbox{Ker}\,\left( \oplus_i L(E_i^R) \longrightarrow
L(E'^R)\right)}{\displaystyle
\mbox{Im}\, \left( \oplus_{\langle a^i,a^j\rangle<\sigma} L(E_i^R\cap E_j^R)
\rightarrow \oplus_i L(E_i^R)\right)} \right)^\ast
\otimes_{\Z}\C \subseteq T^2_Y(-R)
\quad \mbox{(cf. \zitat{3}{3}(2)).}
\]
Then, the $P$-linear map
\bean
^{\displaystyle \kR}\!\!/\!_{\displaystyle \kR_0} &\longrightarrow & A\\
r(a,b;c) & \mapsto &
\left\{ \begin{array}{ll}
\psi_i(a-b)\, x^{\pi(a+c)-R} & \mbox{for } \pi(a)\in K_i^R;\; \pi(a+c)\kgeq R\\
0 & \mbox{for }\pi(a)\kgeq R \mbox{ or } \pi(a+c)\in\bigcup_i K_i^R
\end{array}\right.\\
s(a,b,c) &\mapsto & 0
\eean
is correct defined, and via the formula \zitat{2}{2} it induces the same
element of
$T^2_Y$.
}
\vspace{2ex}
\\
\par
For the prove we refer to \S \ref{s4}. Nevertheless, we check the {\em
correctness of
the definition} of the $P$-linear map
$^{\displaystyle \kR}\!/\!_{\displaystyle \kR_0} \rightarrow A$ instantly:
\begin{itemize}
\item[(i)]
If $\pi(a)$ is contained in two different sets $K_i^R$ and $K_j^R$, then the
two
fundamental generators $a^i$ and $a^j$ can be connected by a sequence
$a^{i_0},\dots,a^{i_p}$, such that
\begin{itemize}
\item[$\bullet$]
$a^{i_0}=a^i,\, a^{i_p}=a^j,$
\item[$\bullet$]
$a^{i_{v-1}}$ and $a^{i_v}$ are the edges of some 2-face of $\sigma$
($v=1,\dots,p$),
and
\item[$\bullet$]
$\pi(a)\in K^R_{i_v}$ for $v=0,\dots,p$.
\end{itemize}
Hence, $\mbox{supp}\,(a-b)\subseteq E^R_{i_{v-1}}\cap E^R_{i_v}$
($v=1,\dots,p$),
and we obtain
\[
\psi_i(a-b)=\psi_{i_1}(a-b)=\dots=\psi_{i_{p-1}}(a-b)=\psi_j(a-b)\,.
\]
\item[(ii)]
There are three types of $P$-linear relations between the generators $r(\dots)$
and
$s(\dots)$ of $\kR$:
\bean
0 &=& \ku{z}^d\,r(a,b;c) -r(a,b;c+d) + r(a+c,b+c;d)\,,\\
0 &=& r(b,c;d) - r(a,c;d) + r(a,b;d) - s(a+d,b+d,c+d) + \ku{z}^d\,
s(a,b,c)\,,\\
0 &=& s(b,c,d) - s(a,c,d) + s(a,b,d) - s(a,b,c)\,.
\eean
Our map respects them all.
\item[(iii)]
Finally, the typical element
$(\ku{z}^a-\ku{z}^b)e^{cd} - (\ku{z}^c-\ku{z}^d)e^{ab} \in \kR_0$
equals
\[
-r(c,d;a)+r(c,d;b)+r(a,b;c)-r(a,b;d) - s(a+c,b+c,a+d) - s(a+d,b+c,b+d)\,.
\]
It will be sent to 0, too.
\vspace{2ex}
\end{itemize}
\par


\neu{36}
The cup product $T^1_Y\times T^1_Y\rightarrow T^2_Y$ respects the grading, i.e.
it splits
into pieces
\[
T^1_Y(-R)\times T^1_Y(-S) \longrightarrow T^2_Y(-R-S)\quad (R,S\in M)\,.
\]
To describe these maps in our combinatorial language, we choose some
set-theoretical
section $\Phi:M\rightarrow\Z^{w+1}$ of the $\Z$-linear map
\bean
\pi: \Z^{w+1} &\longrightarrow& M\\
a&\mapsto&\sum_v a_v\,r^v
\eean
with the additional property $\Phi(\sigma^{\veee}\cap M)\subseteq \N^{w+1}$.\\
\par
Let $q\in L(E)\subseteq\Z^{w+1}$ be an integral relation between the generators
of
the semigroup $\sigma^{\veee}\cap M$. We introduce the following notations:
\begin{itemize}
\item
$q^+,q^-\in\N^{w+1}$ denote the positive and the negative part of $q$,
respectively.
(With other words: $q=q^+-q^-$ and $\sum_v q^-_v\,q^+_v=0$.)
\item
$\bar{q}:=\pi(q^+)=\sum_v q_v^+\,r^v = \pi(q^-)=\sum_v q_v^-\,r^v \in M$.
\item
If $\varphi,\psi: L(E)\rightarrow\Z$ are linear maps and $R,S\in M$, then we
define
\[
t_{\varphi,\psi,R,S}(q):=
\varphi(q)\cdot \psi \left( \Phi(\bar{q}-R)+\Phi(R)-q^-\right) +
\psi(q)\cdot \varphi\left( \Phi(\bar{q}-S)+\Phi(S)-q^+\right)\,.
\vspace{2ex}
\]
\end{itemize}
\par

{\bf Theorem:}
{\em
Assume that $Y$ is smooth in codimension 2.\\
Let $R,S\in M$, and let $\varphi,\psi: L(E)_{\C}\rightarrow\C$ be linear maps
vanishing on $\sum_i L(E_i^R)_{\C}$ and $\sum_i L(E_i^S)_{\C}$, respectively.
In particular,
they define elements $\varphi\in T^1_Y(-R),\,\psi\in T^1_Y(-S)$ (which involves
a slight
abuse of notations).\\
Then, the cup product $\varphi\cup\psi\in T^2_Y(-R-S)$ is given (via
\zitat{3}{3}(3))
by the linear maps $(\varphi\cup\psi)_i: L(E_i^{R+S})_{\C}\rightarrow\C$
defined as follows:
\begin{itemize}
\item[(i)]
If $q\in L(E_i^{R+S})$ (i.e. $\langle a^i,\mbox{supp}\,q\rangle < \langle
a^i,R+S\rangle$) is an integral relation,
then there exists a decomposition $q=\sum_k q^k$ such that
\begin{itemize}
\item
$q^k\in L(E_i^{R+S})$, and moreover
\item
$\langle a^i, \bar{q}^k\rangle < \langle a^i,R+S\rangle$.
\end{itemize}
\item[(ii)]
$(\varphi\cup\psi)_i\left( q\in L(E_i^{R+S})\right):= \sum_k
t_{\varphi,\psi,R,S}(q^k)$.
\vspace{2ex}
\end{itemize}
}
\par

It is even not obvious that the map $q\mapsto \sum_k t(q^k)$
\begin{itemize}
\item
does not depend on the representation of $q$ as a particular sum of $q_k$'s
(which
would instantly imply linearity on $L(E_i^{R+S})$), and
\item
yields the same result on $L(E_i^{R+S}\cap E_j^{R+S})$ for $i,j$ corresponding
to edges
$a^i, a^j$ of some 2-face of $\sigma$.
\end{itemize}
The proof of these facts (cf.\ \zitat{5}{4})and of the entire theorem is
contained
in \S \ref{s5}.\\
\par

{\bf Remark 1:}
Replacing all the terms $\Phi(\bullet)$ in the $t$'s of the previous formula
for
$(\varphi\cup\psi)_i(q)$ by arbitrary liftings from $M$ to $\Z^{w+1}$,
the result in $T^2_Y(-R-S)$ will be unchanged as long as we obey the following
two
rules:
\begin{itemize}
\item[(i)]
Use always (for all $q$, $q^k$, and $i$)
the {\em same liftings} of $R$ and $S$ to $\Z^{w+1}$ (at the places of
$\Phi(R)$ and $\Phi(S)$, respectively).
\item[(ii)]
Elements of $\sigma^{\veee}\cap M$ always have to be lifted to $\N^{w+1}$.
\vspace{2ex}
\end{itemize}

{\bf Proof:}
Replacing $\Phi(R)$ by $\Phi(R)+d$ ($d\in L(E)$) at each occurence changes all
maps $(\varphi\cup\psi)_i$ by the summand $\psi(d)\cdot\varphi(\bullet)$.
However,
this additional linear map comes from $L(E)^\ast$, hence it is trivial on
$\mbox{Ker}\left(\oplus_iL(E_i^{R+S})\rightarrow L(E_0^{R+S})\subseteq
L(E)\right)$.\\
\par
Let us look at the terms $\Phi(\bar{q}-R)$ in $t(q)$ now:
Unless $\bar{q}\kgeq R$, the factor $\varphi(q)$ vanishes (cf. Remark
\zitat{3}{4}). On
the other hand, the expression $t(q)$ is never used for those relations $q$
satisfying
$\bar{q}\kgeq R+S$ (cf. conditions for the $q^k$'s). Hence, we may assume that
\[
(\bar{q}-R)\kgeq 0\; \mbox{ and, moreover, } (\bar{q}-R)\in \bigcup_i K_i^S\,.
\]
Now, each two liftings of $(\bar{q}-R)$ to $\N^{w+1}$ differ by an element of
$\mbox{Ker}\,\psi$ only (apply the method of Remark \zitat{3}{4} again), in
particular,
they cause the same result for $t(q)$.
\hfill$\Box$\\
\par

{\bf Remark 2:}
In the special case of $R\kgeq S\kgeq0$ we can choose liftings $\Phi(R)\geq
\Phi(S)
\geq 0$ in $\N^{w+1}$. Then,
there exists an easier description for $t(q)$:
\begin{itemize}
\item[(i)]
Unless $\bar{q}\kgeq R$, we have $t(q)=0$.
\item[(ii)]
In case of $\bar{q}\kgeq R$ we may assume that $q^+\geq\Phi(R)$ is true
in $\N^{w+1}$.
(The general $q$'s are differences of those ones.) Then, $t$ can be computed as
the
product $t(q)=\varphi(q)\,\psi(q)$.
\vspace{2ex}
\end{itemize}
\par

{\bf Proof:}
(i) As used many times, the property $\bar{q}\in\bigcup_iE_i^R$ implies
$\varphi(q)=0$.
Now, we can distinguish between two cases:\\
{\em Case 1: $\bar{q}\in\bigcup_iE_i^S$.} We obtain $\psi(q)=0$, in particular,
both
summands of $t(q)$ vanish.\\
{\em Case 2: $\bar{q}\kgeq S$.} Then, $\bar{q}-S,\,S\in \sigma^{\veee}\cap M$,
and $\Phi$ lifts
these elements to $\N^{w+1}$. Now, the condition $\bar{q}\in\bigcup_iE_i^R$
implies
that $\varphi\left( \Phi(\bar{q}-S)+\Phi(S)-q^+\right)=0$.\\
\par
(ii)
We can choose $\Phi(\bar{q}-R):=q^+-\Phi(R)$ and
$\Phi(\bar{q}-S):=q^+-\Phi(S)$. Then,
the claim follows straight forward.
\hfill$\Box$\\
\par

%
%
\sect{Proof of the $T^2$-formula}\label{s4}


\neu{41}
We will use the sheaf $\Omega^1_Y=\Omega^1_{A|\C}$ of K\"ahler
differentials for computing the modules $T^i_Y$. The maps
\[
\alpha_i: \mbox{Ext}^i_A\left(
\;^{\displaystyle\Omega_Y^1}\!\!\left/\!_{\displaystyle
\mbox{tors}\,(\Omega_Y^1)}\right.  ,
\, A \right)
\hookrightarrow
\mbox{Ext}^i_A\left(  \Omega^1_Y,\,A\right) \cong T^i_Y\quad
(i=1,2)
\]
are injective. Moreover, they are isomorphisms for
\begin{itemize}
\item
$i=1$, since $Y$ is normal, and for
\item
$i=2$, if $Y$ is smooth in codimension 2.
\vspace{2ex}
\end{itemize}
\par


\neu{42}
As in \cite{T1}, we build a special $A$-free resolution (one step further now)
\[
\CE\stackrel{d_E}{\longrightarrow}\CD\stackrel{d_D}{\longrightarrow}
\CC\stackrel{d_C}{\longrightarrow}\CB \stackrel{d_B}{\longrightarrow}
\;^{\displaystyle\Omega_Y^1}\!\!\left/\!_{\displaystyle
\mbox{tors}\,(\Omega_Y^1)}\right.
\rightarrow 0\,.
\vspace{2ex}
\]
With $L^2(E):=L(L(E))$, $L^3(E):=L(L^2(E))$, and
\[
\mbox{supp}^2\xi:= \bigcup_{q\in supp\,\xi} \mbox{supp}\,q\quad (\xi\in
L^2(E)),\quad
\mbox{supp}^3\omega:= \bigcup_{\xi\in supp\,\omega}\mbox{supp}^2\xi\quad
(\omega\in L^3(E)),
\]
the $A$-modules involved in this resolution are defined as follows:
\[
\begin{array}{rcl}
\CB&:=&\oplus_{r\in E} \,A\cdot B(r),\qquad
\CC\,:=\,\oplus_{\!\!\!\!\!\begin{array}[b]{c}\ks q\in L(E) \vspace{-1ex}\\
\ks\ell\kgeq supp\, q\end{array}}
\!\!\!A\cdot C(q;\ell),\\
\CD&:=&\left(
\oplus_{\!\!\!\!\!\!\!\begin{array}{c}\ks q\in L(E)\vspace{-1ex}\\
\ks\eta\kgeq\ell\kgeq supp\, q\end{array}}
\!\!\!\!A\cdot D(q;\ell,\eta) \right)
\oplus \left(
\oplus_{\!\!\!\!\begin{array}{c}\ks\xi\in L^2(E)\vspace{-1ex}\\
\ks\eta\kgeq supp^2 \xi\end{array}}
\!\!\!A\cdot D(\xi;\eta) \right),\;\mbox{ and}\\
\CE&:=&
\begin{array}[t]{r} \left(
\oplus_{\!\!\!\!\!\!\!\!\begin{array}{c}\ks q\in L(E)\vspace{-1ex}\\
\ks\mu\kgeq\eta\kgeq\ell\kgeq supp\, q\end{array}}
\!\!\!\!\!A\cdot E(q;\ell,\eta,\mu) \right)
\oplus \left(
\oplus_{\!\!\!\!\!\!\!\begin{array}{c}\ks\xi\in L^2(E)\vspace{-1ex}\\
\ks\mu\kgeq\eta\kgeq supp^2 \xi\end{array}}
\!\!\!A\cdot E(\xi;\eta,\mu) \right) \oplus \qquad \\
\oplus \left(
\oplus_{\!\!\!\!\begin{array}{c}\ks\omega\in L^3(E)\vspace{-1ex}\\
\ks\omega\kgeq supp^3 \omega\end{array}}
\!\!\!\! A\cdot E(\omega;\mu)\right)
\end{array}
\end{array}
\]
($B,C,D,$ and $E$ are just symbols).
The differentials equal
\[
\begin{array}{cccl}
d_B: &B(r)&\mapsto &d\,x^r\vspace{1ex}\\
d_C: &C(q;\ell)&\mapsto &\sum_{r\in E} q_r\,x^{\ell-r}\cdot B(r)\vspace{1ex}\\
d_D: &D(q;\ell,\eta)&\mapsto &C(q;\eta) - x^{\eta-\ell}\cdot C(q,\ell)\\
d_D: &D(\xi;\eta)&\mapsto& \sum_{q\in L(E)}\xi_q\cdot C(q,\eta)\vspace{1ex}\\
d_E: &E(q;\ell,\eta,\mu)&\mapsto& D(q;\eta,\mu)-D(q;\ell,\mu)+
x^{\mu-\eta}\cdot D(q;\ell,\eta) \\
d_E: &E(\xi;\eta,\mu)&\mapsto &D(\xi;\mu) - x^{\mu-\eta}\cdot D(\xi;\eta) -
\sum_{q\in L(E)} \xi_q\cdot D(q;\eta,\mu)\\
d_E: &E(\omega;\mu)&\mapsto &\sum_{\xi\in L^2(E)} \omega_{\xi}\cdot
D(\xi;\mu)\, .
\vspace{2ex}
\end{array}
\]
\par

Looking at these maps, we see that the complex is $M$-graded: The degree of
each of
the elements $B$, $C$, $D$, or $E$ can be obtained by taking the last of its
parameters
($r$, $\ell$, $\eta$, or $\mu$, respectively).\\
\par

{\bf Remark:} If one prefered a resolution with free $A$-modules of finite rank
(as it was
used in
\cite{T1}), the following replacements would be necessary:
\begin{itemize}
\item[(i)]
Define succesively $F\subseteq L(E)$, $G\subseteq L(F) \subseteq L^2(E)$, and
$H\subseteq L(G)\subseteq L^2(F) \subseteq L^3(E)$ as the finite
sets of normalized, minimal relations between elements of $E$, $F$, or
$G$, respectively. Then, use them instead of $L^i(E)$ ($i=1,2,3$).
\item[(ii)]
Let $\ell$, $\eta$, and $\mu$ run through finite generating
(under $(\sigma^{\veee}\cap M)$-action)
systems of all possible elements meeting the desired inequalities.
\end{itemize}
The disadvantages of those treatment are a more comlplicated description of
the resolution, on the one hand, and difficulties to obtain the
commutative diagram \zitat{4}{3},
on the other hand.\\
\par


\neu{43}
Combining the two exact sequences
\[
^{\kd \kR}\!/\!_{\kd I\,\kR} \longrightarrow
A^m \longrightarrow
^{\kd I}\!\!/\!_{\kd I^2}\rightarrow 0\quad
\mbox{and}\quad
^{\kd I}\!\!/\!_{\kd I^2}\longrightarrow
\Omega^1_{\C^{w+1}}\otimes A \longrightarrow
\Omega_Y^1 \rightarrow 0\,,
\]
we get a complex (not exact at the place of $A^m$) involving $\Omega_Y^1$. We
will compare
in the following commutative diagram this complex with the previous resolution
of
$^{\displaystyle\Omega_Y^1}\!\!\left/\!_{\displaystyle
\mbox{tors}\,(\Omega_Y^1)}\right.$:
\vspace{-5ex}\\
\[
\dgARROWLENGTH=0.8em
\begin{diagram}
\node[5]{^{\kd I}\!\!/\!_{\kd I^2}}
\arrow{se,t}{d}\\
\node[2]{^{\kd \kR}\!/\!_{\kd I\,\kR}}
\arrow[2]{e}
\arrow{se,t}{p_D}
\node[2]{A^m}
\arrow{ne}
\arrow[2]{e}
\arrow[2]{s,l}{p_C}
\node[2]{\Omega_{\C^{w+1}}\!\otimes \!A}
\arrow{e}
\arrow[2]{s,lr}{p_B}{\sim}
\node{\Omega_Y}
\arrow[2]{s}
\arrow{e}
\node{0}\\
\node[3]{\mbox{Im}\,d_D}
\arrow{se}\\
\node{\CE}
\arrow{e,t}{d_E}
\node{\CD}
\arrow[2]{e,t}{d_D}
\arrow{ne}
\node[2]{\CC}
\arrow[2]{e,t}{d_C}
\node[2]{\CB}
\arrow{e}
\node{^{\displaystyle\Omega_Y^1}\!\!\!\left/\!\!_{\displaystyle
\mbox{tors}\,(\Omega_Y^1)}\right.}
\arrow{e}
\node{0}
\end{diagram}
\]
\par

Let us explain the three labeled vertical maps:
\begin{itemize}
\item[(B)]
$p_B: dz_r \mapsto B(r)$ is an isomorphism between two free $A$-modules of rank
$w+1$.
\item[(C)]
$p_C: e^{ab} \mapsto C(a-b;\pi(a))$. In particular, the image of this map is
spanned by
those $C(q,\ell)$ meeting $\ell\kgeq \bar{q}$ (which is stronger than just
$\ell\kgeq\mbox{supp}\,q$).
\item[(D)]
Finally, $p_D$ arises as pull back of $p_C$ to $^{\kd \kR}\!/\!_{\kd I\,\kR}$.
It can
be described by
$r(a,b;c)\mapsto D(a-b; \pi(a),\pi(a+c))$ and $s(a,b,c)\mapsto D(\xi;\pi(a))$
($\xi$ denotes
the relation $\xi=[(b-c)-(a-c)+(a-b)=0]$).
\vspace{2ex}
\end{itemize}
\par

{\bf Remark:}
Starting with the typical $\kR_0$-element mentioned in \zitat{3}{5}(iii), the
previous
description of the map $p_D$ yields 0 (even in $\CD$).\\
\par


\neu{44}
By \zitat{4}{1} we get the $A$-modules $T^i_Y$ by computing the cohomology
of the complex dual to those of \zitat{4}{2}.\\
\par
As in \cite{T1}, denote by $G$ one of the capital letters $B$, $C$, $D$, or
$E$. Then, an element $\psi$ of the dual free module $(\bigoplus\limits_G
\C[\check{\sigma}\cap M]\cdot G)^\ast$ can be described by giving elements
$g(x)\in\C[\check{\sigma}\cap M]$ to be the images of the generators $G$
($g$ stands for $b$, $c$, $d$, or $e$, respectively).\\
\par
For $\psi$ to be homogeneous of degree $-R\in M$, $g(x)$ has to be
a monomial of degree
\[
\deg g(x)=-R+\deg G.
\]
In particular, the corresponding complex coefficient $g\in \C$ (i.e.
$g(x)=g\cdot x^{-R+\deg G}$) admits the property that
\[
g\neq 0\quad\mbox{implies}\quad -R+\deg G\ge 0\quad (\mbox{i.e.}\;
-R+\deg G\in\check{\sigma}).
\vspace{2ex}
\]
\par


{\bf Remark:}
Using these notations,
Theorem \zitat{3}{3}(1) was proved in \cite{T1} by showing that
\bean
\left(\left. ^{\displaystyle L(E_0^R)}\!\!\right/
\!_{\displaystyle \sum_i L(E_i^R)} \right)^\ast \otimes_{\Z} \C
&\longrightarrow&
^{\kd \mbox{Ker}(\CC^\ast_{-R}\rightarrow \CD^\ast_{-R})}\!\!\left/
\!_{\kd \mbox{Im}(\CB^\ast_{-R}\rightarrow\CC^\ast_{-R})}\right.\\
\varphi &\mapsto&
[\dots,\, c(q;\ell):=\varphi(q),\dots]
\eean
is an isomorphism of vector spaces.\\
Moreover, looking at the diagram of
\zitat{4}{3}, $e^{ab}\in A^m$ maps to both $\ku{z}^a-\ku{z}^b\in ^{\kd
I}\!\!/\!
_{\kd I^2}$ and $C(a-b;\pi(a))\in \CC$. In particular, we can verify Theorem
\zitat{3}{4}:
Each $\varphi:L(E)_{\C} \rightarrow\C$, on the one hand,
and its associated $A$-linear map
\bean
^{\displaystyle I}\!\!/\!_{\displaystyle I^2} &\longrightarrow& A\\
\ku{z}^a-\ku{z}^b & \mapsto & \varphi(a-b)\cdot x^{\pi(a)-R},
\eean
on the other hand, induce the same element of $T^1_Y(-R)$.  \\
\par


\neu{45}
For computing $T_Y^2(-R)$,
the interesting part of the dualized complex
$\zitat{4}{2}^\ast$ in degree $-R$ equals the complex
of $\C$-vector spaces
\[
\CC^{\ast}_{-R} \stackrel{d_D^{\ast}}{\longrightarrow} \CD^{\ast}_{-R}
\stackrel{d_E^{\ast}}{\longrightarrow} \CE^{\ast}_{-R}
\]
with coordinates $\underline{c}$, $\underline{d}$, and $\underline{e}$,
respectively:
\bean
\CC^{\ast}_{-R} &=&
\{\underline{c(q;\ell)}\, |\;
c(q;\ell)=0 \;\mbox{for}\;\ell-R\notin\check{\sigma}\}\\
\CD^{\ast}_{-R} &=&
\{[\underline{d(q;\ell,\eta)},\underline{d(\xi;\eta)}]\;|\;
\begin{array}[t]{ccccl}
d(q;\ell,\eta)&=&0& \mbox{for} &\eta-R\notin\check{\sigma}\mbox{, and}\\
d(\xi;\eta)&=&0& \mbox{for} &\eta-R\notin\check{\sigma} \}
\end{array}\\
\CE^{\ast}_{-R} &=&
\{ [\underline{e(q;\ell,\eta,\mu)},
\underline{e(\xi;\eta,\mu)},
\underline{e(\omega ;\mu)}]\,|\;
\mbox{each coordinate vanishes for } \mu - R \notin \check{\sigma} \}.
\vspace{1ex}
\eean
\par

The differentials
$d_D^{\ast}$ and $d_E^{\ast}$ can be described by
\[
\begin{array}{lcll}
d(q;\ell,\eta)&=&c(q;\eta)-c(q;\ell),&\\
d(\xi;\eta)&=&\sum\limits_{q\in F}\xi_q\cdot c(q;\eta),&
\mbox{and}\\
e(q;\ell,\eta,\mu) &=&
d(q;\eta,\mu) - d(q;\ell,\mu) + d(q;\ell,\eta),\\
e(\xi;\eta,\mu) &=&
d(\xi;\mu) - d(\xi;\eta) - \sum_{q\in F} \xi_q\cdot
d(q;\eta,\mu),\\
e(\omega ;\mu) &=&
\sum\limits_{\xi\in G} \omega_{\xi}\cdot d(\xi;\mu).
\end{array}
\vspace{1ex}
\]
\par

Denote $V:= \mbox{Ker}\,d^{\ast}_E \subseteq \CD_{-R}^{\ast}\,$ and
$\,W:= \mbox{Im}\,d_D^{\ast}\subseteq V$, i.e.
\bean
V&=& \{ [\underline{d(q;\ell,\eta)};\,\underline{d(\xi;\eta)}]\,|\;
\begin{array}[t]{l}
q\in L(E), \;\eta\kgeq\ell\kgeq\mbox{supp}\,q\mbox{ in }M;\\
\xi\in L^2(E), \;\eta\kgeq\mbox{supp}^2\xi;
\vspace{0.5ex}\\
d(q;\ell,\eta) = d(\xi;\eta) = 0 \mbox{ for } \eta -R \notin \check{\sigma},\\
d(q;\ell,\mu) = d(q;\ell,\eta) + d(q;\eta,\mu) \; (\mu\kgeq\eta\kgeq\ell\kgeq
\supp q),\\
d(\xi;\mu)= d(\xi;\eta) + \sum_q \xi_q \cdot d(q;\eta,\mu)\;
(\mu\kgeq\eta\kgeq \mbox{supp}^2 \xi),\\
\sum_{\xi\in G}\omega_{\xi} \,d(\xi;\mu) =0 \mbox{ for }\omega \in L^3(E)
\mbox{ with } \mu \kgeq \mbox{supp}^3\omega\,\},
\end{array}\\
W&=& \{ [\underline{d(q;\ell,\eta)};\,\underline{d(\xi;\eta)}]\,|\;
\mbox{there are $c(q;\ell)$'s with}
\begin{array}[t]{l}
c(q,\ell)=0 \mbox{ for } \ell-R\notin\check{\sigma},\\
d(q;\ell,\eta) = c(q;\eta)-c(q;\ell),\\
d(\xi;\eta)= \sum_q\xi_q\cdot c(q;\eta)\,\}.
\end{array}
\eean
By construction, we obtain
\[
V\!\left/\!_{\displaystyle W}\right. =
\mbox{Ext}^i_A\left(
\;^{\displaystyle\Omega_Y^1}\!\!\left/\!_{\displaystyle
\mbox{tors}\,(\Omega_Y^1)}\right.
, \, A \right)(-R)
\subseteq T^2_Y(-R)
\]
(which is an isomorphism, if $Y$ is smooth in codimension 2).\\
\par


\neu{46}
Let us define the much easier vector spaces
\bean
V_1&:=& \{[\underline{x_i(q)}_{(q\in L(E_i^R))}]\,|\;
\begin{array}[t]{l}
x_i(q)=x_j(q) \mbox{ for }
\begin{array}[t]{l}
\bullet\,
\langle a^i, a^j \rangle < \sigma \mbox{ is a 2-face and}\\
\bullet\,
q\in L(E_i^R\cap E_j^R)\,,
\end{array}\\
\xi\in L^2(E_i^R) \mbox{ implies } \sum_q\xi_q\cdot
x_i(q)=0
\,\}\;\mbox{ and}
\end{array}
\vspace{1ex}
\\
W_1&:=& \{[\underline{x(q)}_{(q\in \cup_i L(E_i^R))}]
\,|\;
\begin{array}[t]{l}
\xi\in L(\bigcup_i L(E_i^R)) \mbox{ implies } \sum_q\xi_q\cdot
x(q)=0 \,\}.
\vspace{2ex}
\end{array}
\eean
\par

{\bf Lemma:}
{\em
The linear map $V_1\rightarrow V$ defined by
\bean
d(q;\ell,\eta) &:=& \left\{
\begin{array}{lll}
x_i(q) &\mbox{ for }& \ell\in K_i^R,\;\; \eta \kgeq R\\
0 &\mbox{ for }& \ell \kgeq R \;\mbox{ or } \;\eta \in
\bigcup_i K_i^R\,;
\end{array} \right.  \\
d(\xi;\eta) &:=&
0\,
\eean
induces an injective map
\[
V_1\!\left/\!_{\displaystyle W_1}\right.
\hookrightarrow
V\!\left/\!_{\displaystyle W}\right.\,.
\]
If $Y$ is smooth in codimension 2, it will be an isomorphism.
}
\\
\par

{\bf Proof:}
1) The map $V_1 \rightarrow V$ is {\em correct defined}:
On the one hand, an argument as used in \zitat{3}{5}(i) shows that $\ell\in
K_i^R\cap K_j^R$ would imply $x_i(q)=x_j(q)$. On the other hand,
the image of
$[x_i(q)]_{q\in L(E_i^R)}$ meets all conditions in the definition of $V$.
\vspace{1ex}
\\
2) $W_1$ maps to $W$ (take $c(q,\ell):=x(q)$ for $\ell\kgeq R$ and
$c(q,\ell):=0$ otherwise).
\vspace{1ex}
\\
3) The map between the two factor spaces is {\em injective}: Assume for
$[x_i(q)]_{q\in L(E_i^R)}$ that there exist elements $c(q,\ell)$, such that
\bean
c(q;\ell) &=& 0 \; \mbox{ for } \ell \in \bigcup_i K^R_i\, ,\\
x_i(q) &=& c(q;\eta) - c(q;\ell) \; \mbox{ for }
\eta \kgeq \ell,\, \ell\in K_i^R,\, \eta\kgeq R\,,\\
0 &=&
c(q;\eta) - c(q;\ell)\; \mbox{ for } \eta\kgeq \ell \mbox{ and }
[\ell\kgeq R\mbox{ or } \eta\in \bigcup_iK_i^R]\, , \mbox{ and}\\
0 &=&
\sum_q \xi_q \cdot c(q;\eta)  \; \mbox{ for } \eta \kgeq
\mbox{supp}^2\xi\, .
\eean
In particular, $x_i(q)$ do not depend on $i$, and these elements
meet the property
\[
\sum_q \xi_q \cdot x_{\bullet}(q) = 0 \; \mbox{ for } \xi\in L(\bigcup_i
L(E_i^R)).
\]
4) If $Y$ is smooth in codimension 2, the map is {\em surjective} :\\
Given an element $[d(q;\ell,\eta),\,d(\xi;\eta)]\in V$, there exist
complex numbers $c(q;\eta)$ such that:
\begin{itemize}
\item[(i)]
$d(\xi;\eta) = \sum_q\xi_q\cdot c(q;\eta)\,$ ,
\item[(ii)]
$c(q;\eta)=0 \mbox{ for } \eta\notin R+\sigma^{\veee}\,
(\mbox{i.e. }\eta\in \bigcup_iK_i^R)\,$.
\end{itemize}
(Do this separately for each $\eta$ and distinguish between the cases
$\eta\in R +\sigma^{\veee}$ and $\eta\notin R+\sigma^{\veee}$.)\\
In particular, $[c(q;\eta) - c(q;\ell),\, d(\xi;\eta)]\in W$. Hence, we
have seen that we may assume $d(\xi;\eta)=0$.\\
\par
Let us choose some sufficiently high degree $\ell^\ast\kgeq E$.
Then,
\[
x_i(q):= d(q;\ell,\eta) - d(q;\ell^\ast\!,\eta)
\]
(with $\ell\in K_i^R$, $\ell\kgeq \mbox{supp}\,q$
(cf.\ Lemma \zitat{3}{2}(2)), and $\eta\kgeq\ell,\ell^\ast\!,R$)
defines some preimage:
\begin{itemize}
\item[(i)]
It is independent from the choice of $\eta$: Using a different $\eta'$
generates
the difference $d(q;\eta,\eta')-d(q;\eta,\eta')$.
\item[(ii)]
It is independent from $\ell\in K_i^R$: Choosing another $\ell'\in K_i^R$
with $\ell'\kgeq\ell$ would add the summand $d(q;\ell,\ell')$, which is 0;
for the general case use Lemma \zitat{3}{2}(2).
\item[(iii)]
If $\langle a^i,a^j\rangle < \sigma$ is a 2-face with $\mbox{supp}\,q
\subseteq L(E^R_i)\cap L(E_j^R)$, then by Lemma \zitat{3}{2}(2) we can choose
an
$\ell\in K_i^R\cap K_j^R$ achieving $x_i(q)=x_j(q)$.
\item[(iv)]
For $\xi\in L^2(E_i^R)$ we have
\[
\sum_q \xi_q\cdot d(q;\ell,\eta) = \sum_q \xi_q\cdot d(q;\ell^\ast\!,\eta) =
0\,,
\]
and this gives the corresponding relation for the $x_i(q)'$s.
\item[(v)]
Finally, if we apply to
$[\ku{x_i(q)}]\in V_1$
the linear map $V_1\rightarrow V$, the result differs from
$[d(q;\ell,\eta),0]\in V$ by
the $W$-element built from
\[
c(q;\ell) := \left\{ \begin{array}{ll}
d(q;\ell,\eta) - d(q;\ell^\ast\!,\eta) & \mbox{ if } \ell\kgeq R \\
0 & \mbox{ otherwise }.
\end{array} \right.
\vspace{-2ex}
\]
\end{itemize}
\hfill$\Box$\\
\par


\neu{47}
Now, it is easy to complete the proofs for Theorem \zitat{3}{3} (part 2 and 3)
and
Theorem \zitat{3}{5}:\\
\par
First, for a tuple $[\ku{x_i(q)}]_{q\in L(E_i^R)}$, the condition
\[
\xi\in L^2(E_i^R) \mbox{ implies } \sum_q\xi_q\cdot x_i(q)=0
\]
is equivalent to the fact the components $x_i(q)$ are induced by elements
$x_i\in L(E_i^R)_{\C}^\ast$.\\
The other condition for elements of $V_1$ just says that for 2-faces
$\langle a^i,a^j\rangle<\sigma$ there is $x_i=x_j$ on
$L(E_i^R\cap E_j^R)_{\C}=L(E_i^R)_{\C}\cap L(E_j^R)_{\C}$. In particular, we
obtain
\[
V_1= \mbox{Ker}\left( \oplus_i L(E_i^R)_{\C}^\ast \rightarrow
\oplus_{\langle a^i,a^j\rangle <\sigma} L(E_i^R\cap E_j^R)_{\C}^\ast \right)\,.
\]
In the same way we get
\[
W_1 = \left( \sum_i L(E^R_i)_{\C}\right)^\ast\,,
\]
and our $T^2$-formula is proven.\\
\par
Finally, if $\psi_i:L(E_i^R)_{\C}\rightarrow \C$ are linear maps defining an
element of
$V_1$, they induce the following $A$-linear map on $\CD$ (even on
$\mbox{Im}\,d_D$):
\bean
D(q;\ell,\eta) &\mapsto& \left\{
\begin{array}{lll}
\psi_i(q)\cdot x^{\eta-R} &\mbox{ for }& \ell\in K_i^R,\;\; \eta \kgeq R\\
0 &\mbox{ for }& \ell \kgeq R \;\mbox{ or } \;\eta \in
\bigcup_i K_i^R
\end{array} \right.\\
D(\xi;\eta &\mapsto& 0\,.
\eean
Now, looking at the diagram of \zitat{4}{3}, this translates exactly into the
claim of
Theorem \zitat{3}{5}.\\
\par

%
%
\sect{Proof of the cup product formula}\label{s5}


\neu{51}
Fix an $R\in M$, and let $\varphi\in L(E)^\ast_{\C}$ induce some element
(also denoted by $\varphi$)
of $T^1_Y(-R)$. Using the notations of \zitat{2}{3}, \zitat{3}{4},
and \zitat{3}{6} we can take
\[
\tilde{\varphi}(f_{\alpha\beta}):=
\varphi(\alpha-\beta)\cdot \ku{z}^{\Phi(\pi(\alpha)-R)}
\]
for the auxiliary $P$-elements needed to compute the $\lambda(\varphi)$'s
(cf. Theorem \zitat{3}{4}).\\
\par
Now, we have to distinguish between the two several types of relations
generating the $P$-module $\kR\subseteq P^m$:
\begin{itemize}
\item[(r)]
Regarding the relation $r(a,b;c)$ we obtain
\bean
\sum_{(\alpha,\beta)\in m} r(a,b;c)_{\alpha\beta}\cdot
\tilde{\varphi}(f_{\alpha\beta}) &=&
\tilde{\varphi}(f_{a+c,b+c}) - \ku{z}^c\,\tilde{\varphi}(f_{ab})
\\
&=&
\varphi(a-b)\cdot \left(
\ku{z}^{\Phi(\pi(a+c)-R)} - \ku{z}^{c+\Phi(\pi(a)-R)} \right)
\\
&=&
\varphi(a-b)\cdot
f_{\Phi(\pi(a+c)-R),\,c+\Phi(\pi(a)-R)}\,.
\eean
In particular,
\[
\lambda_{\alpha\beta}^{r(a,b;c)}(\varphi) =
\left\{\begin{array}{ll}
\varphi(a-b) & \mbox{ for } [\alpha,\beta] = [c+\Phi(\pi(a)-R),\,
\Phi(\pi(a+c)-R)]\\
0 & \mbox{ otherwise}\,.
\end{array} \right.
\]
\item[(s)]
The corresponding result for the relation $s(a,b,c)$ is much nicer:
\bean
\sum_{(\alpha,\beta)\in m} s(a,b,c)_{\alpha\beta}\cdot
\tilde{\varphi}(f_{\alpha\beta}) &=&
\tilde{\varphi}(f_{bc})-
\tilde{\varphi}(f_{ac})+
\tilde{\varphi}(f_{ab})\\
&=&
[\varphi(b-c)-\varphi(a-c)+\varphi(a-b)]\cdot
\ku{z}^{\Phi(\pi(a)-R)}\\
&=& 0\,.
\eean
In particular, $\lambda^{s(a,b,c)}(\varphi)=0$.
\vspace{2ex}
\end{itemize}
\par


\neu{52}
Now, let $R,S,\varphi$, and $\psi$ as in the assumption of Theorem
\zitat{3}{6}. Using formula \zitat{2}{3}(iii), our previous computations
yield $(\varphi\cup\psi)(s(a,b,c))=0$ and
\[
\begin{array}{l}
(\varphi\cup\psi)(r(a,b;c))=
\sum_{\alpha,\beta}\lambda^{r(a,b;c)}_{\alpha\beta}(\varphi)
\cdot \psi(f_{\alpha\beta}) +
\sum_{\alpha,\beta} \varphi(f_{\alpha\beta})\cdot
\lambda^{r(a,b;c)}_{\alpha\beta}(\psi)
\vspace{2ex}\\
\qquad=
\begin{array}[t]{r}
\varphi(a-b)\cdot \psi\left(
c^{}+\Phi(\pi(a)-R)-\Phi(\pi(a+c)-R)\right)\cdot
x^{\pi(c+\Phi(\pi(a)-R))-S} +\qquad\\
+\psi(a-b)\cdot \varphi\left(
c+\Phi(\pi(a)-S)-\Phi(\pi(a+c)-S)\right)\cdot
x^{\pi(c+\Phi(\pi(a)-S))-R}
\end{array}
\vspace{2ex}\\
\qquad=
\begin{array}[t]{r}
\left[ \varphi(a-b)\cdot \psi\left(
c+\Phi(\pi(a)-R)-\Phi(\pi(a+c)-R)\right) +\right.
\qquad\qquad\qquad\qquad\qquad\\
\left. + \psi(a-b)\cdot \varphi\left(
c+\Phi(\pi(a)-S)-\Phi(\pi(a+c)-S)\right)
\right]
\cdot x^{\pi(a+c)-R-S}\,.
\vspace{1ex}
\end{array}
\end{array}
\]
\par

{\bf Remark:}
Unless $\pi(a+c)\kgeq R+S$, both summand in the brackets will vanish. For
instance,
on the one hand, $\pi(a)\in\bigcup_iK_i^R$ would cause $\varphi(a-b)=0$, and,
on the
other hand, $\pi(a)-R\kgeq 0$ and $\pi(c+\Phi(\pi(a)-R))\in \bigcup_iK_i^S$
imply $\psi(c+\Phi(\pi(a)-R)-\Phi(\pi(a+c)-R))=0$.\\
\par

To apply Theorem \zitat{3}{5} we would like to remove the argument $c$ from
the big coefficient. This will be done by adding a suitable
coboundary $T$.\\
\par


\neu{53}
Let us start with defining for $(\alpha,\beta)\in m$
\[
t(\alpha,\beta):= \begin{array}[t]{r}
\varphi(\alpha-\beta)\cdot
\psi \left( \Phi(\pi(\alpha)-R)+\Phi(R)-\beta\right)+\qquad\qquad\qquad\\
+ \psi(\alpha-\beta) \cdot
\varphi \left( \Phi(\pi(\alpha)-S)+\Phi(S)-\alpha\right)\,.
\end{array}
\]
(This expression is related to $t_{\varphi,\psi,R,S}$ from
\zitat{3}{6} by $t(q)=t(q^+,q^-)$.) \\
\par

{\bf Lemma:}{\em
Let $\alpha,\beta,\gamma\in\N^E$ with $\pi(\alpha)=\pi(\beta)=\pi(\gamma)$.
\begin{itemize}
\item[(1)]
$t(\alpha,\beta)=t(\alpha-\beta)$
as long as $\pi(\alpha)\in\bigcup_i K_i^{R+S}$.
\item[(2)]
$t(\beta,\gamma)-t(\alpha,\gamma)+t(\alpha,\beta)=0$.
\vspace{2ex}
\end{itemize}
}
\par

{\bf Proof:}
(1) It is enough to show that $t(\alpha+r,\beta+r)=t(\alpha,\beta)$ for
$r\in \N^E$, $\pi(\alpha+r)\in\bigcup_iK_i^{R+S}$. But the difference of these
two terms has exactly the shape of the coefficient of $x^{\pi(a+c)-R-S}$ in
\zitat{5}{2}. In particular, the argument given in the previous remark
applies again.\\
\par
(2) By extending $\varphi$ and $\psi$ to linear maps $\C^E\rightarrow\C$,
we obtain
\[
t(\alpha,\beta) = \begin{array}[t]{r}
[\varphi(\alpha-\beta)\,\psi\left(
\Phi(\pi(\alpha)-R) + \Phi(R)\right) + \psi(\alpha-\beta) \,
\varphi\left( \Phi(\pi(\alpha)-S)+\Phi(S)\right)]+\,\\
+[\varphi(\beta)\,\psi(\beta)-\varphi(\alpha)\,\psi(\alpha)].
\end{array}
\]
Now, since $\pi(\alpha)=\pi(\beta)=\pi(\gamma)$, both types of summands add
up to 0 separately in
$t(\beta,\gamma)-t(\alpha,\gamma)+t(\alpha,\beta)$.
\hfill$\Box$\\
\par

{\bf Remark:} The previous lemma does not imply that $t(q)$ is
$\Z$-linear in $q$. The
assumption for $\pi(\alpha)$ made in (1) is really essential.\\
\par

Now, we obtain a
$P$-linear map $T\in \mbox{Hom}(P^m,A)$ by
\[
T: e^{\alpha\beta}\mapsto
\left\{ \begin{array}{ll}
t(\alpha,\beta)\,x^{\pi(\alpha)-R-S} & \mbox{ for } \pi(\alpha)\kgeq R+S\\
0 & \mbox{ otherwise}\,.
\end{array} \right.
\]
Pulling back $T$ to $\kR\subseteq P^m$ yields (in case of $\pi(a+c)\kgeq R+S$)
\bean
T(r(a,b;c)) &=& \left\{ \begin{array}{ll}
[t(a+c,b+c)-t(a,b)]\cdot x^{\pi(a+c)-R-S} & \mbox{ for } \pi(a)\kgeq R+S\\
t(a+c,b+c)\cdot x^{\pi(a+c)-R-S} & \mbox{ otherwise}
\end{array} \right.\\
&=&
\left\{ \begin{array}{ll}
-(\varphi\cup\psi)(r(a,b;c)) & \mbox{ for } \pi(a)\kgeq R+S\\
t(a,b)\,x^{\pi(a+c)-R-S} -(\varphi\cup\psi)(r(a,b;c)) & \mbox{ otherwise}\,
\end{array} \right.
\eean
and $T(s(a,b,c))=0$ (by (2) of the previous lemma).\\
\par

On the other hand, $T$ yields a trivial element of $T^2_Y(-R-S)$,
i.e. inside this group we may replace
$\varphi\cup\psi$ by $(\varphi\cup\psi)+T$ to obtain
\bean
(\varphi\cup\psi)(r(a,b;c)) &=&
\left\{ \begin{array}{ll}
t(a,b)\cdot x^{\pi(a+c)-R-S}&\mbox{ for } \pi(a)\in \bigcup_i K_i^{R+S};\;
\pi(a+c)\kgeq R+S\\
0 & \mbox{ otherwise}\,,
\end{array} \right.
\vspace{1ex}\\
(\varphi\cup\psi)(s(a,b,c)) &=& 0\,.
\vspace{1ex}
\eean
\par

Having Theorem \zitat{3}{5} in mind, this formula for $\varphi\cup\psi$ is
exactly what
we were looking for:\\
Given an $r(a,b;c)$ with $\pi(a)\in K_i^{R+S}$,
let us compute $(\varphi\cup\psi)_i(q:=a-b)$ following the recipe of (i), (ii)
of Theorem
\zitat{3}{6}. We do not need to split
$q=a-b$ into a sum $q=\sum_k q^k$ - the element $q$ itself already satisfies
the condition
\[
\langle a^i,\bar{q}\rangle \leq \langle a^i, \pi(a) \rangle < \langle
a^i,R+S\rangle.
\]
In particular, with $(\varphi\cup\psi)_i(a-b)=t(a-b)=t(a,b)$ we
will obtain the right result - if the recipe is assumed to be correct. \\
\par


\neu{54}
We will fill those remaining gap now, i.e. we will show that
\begin{itemize}
\item[(a)]
each $q\in L(E_i^{R+S})$
admits a decomposition $q=\sum_k q^k$ with the desired properties,
\item[(b)]
$\sum_k q^k=0$ (with $\bar{q}^k\in K_i^{R+S}$) implies $\sum_k t(q^k)=0$, and
\item[(c)]
for adjacent $a^i,a^j$ the relations $q\in L(E_i^{R+S}\cap E_j^{R+S})$ admit
a decomposition $q=\sum_kq^k$ that works for both $i$ and $j$.
\end{itemize}
(In particular, this answers the questions arised right after stating the
theorem in
\zitat{3}{6}.)\\
\par

Let us fix an element $i\in \{1,\dots,N\}$. Since $\sigma^{\veee}\cap M$
contains elements $r$ with $\langle a^i,r\rangle =1$, some of them must be
contained in the generating set $E$, too. We choose one of these elements
and call it $r(i)$.\\
Now, to each $r\in E$ we associate some relation $p(r)\in L(E)$ by
\[
p(r):= e^r - \langle a^i, r \rangle\cdot e^{r(i)} +
[\mbox{suitable element of } \Z^{E\cap (a^i)^\bot}]\,.
\]
The two essential properties of these special relations are
\begin{itemize}
\item[(i)]
$\langle a^i, \bar{p}(r)\rangle = \langle a^i, r\rangle$, and
\item[(ii)]
if $q\in L(E)$ is any relation, then $q$ and $\sum_{r\in E}q_r\cdot p(r)$
differ
by some element of $L(E\cap (a^i)^\bot)$ only.
\vspace{1ex}
\end{itemize}
\par

In particular, this proves (a). For (b) we start with the following\\
\par

{\em Claim:}
Let $q^k\in L(E)$ be relations such that
$\sum_k \langle a^i,\bar{q}^k\rangle < \langle a^i, R+S\rangle$.
Then, $\sum_k t(q^k)=t(\sum_k q^k)$.\\
\par
{\em Proof:} We can restrict ourselves to the case of two summands, $q^1$ and
$q^2$. Then,
by Lemma \zitat{5}{3},
\bean
t(q^1)+t(q^2) &=&
t\left((q^1)^+,(q^1)^-\right) + t\left((q^2)^+,(q^2)^-\right)\\
&=&
t\left((q^1)^++(q^2)^+,(q^1)^-+(q^2)^+\right) +
t\left((q^2)^++(q^1)^-,(q^2)^-+(q^1)^-\right)\\
&=&
t\left((q^1)^++(q^2)^+,(q^2)^-+(q^1)^-\right)\\
&=& t(q^1+q^2)\,.
\hspace{9cm} \Box
\eean
\par

In particular, if $\sum_kq^k=0$ (with $\bar{q}^k\in K_i^{R+S}$), then this
applies for
the special decompositions
\[
q^k=\sum_r q^k_r\cdot p(r) + q^{0,k} \quad (q^{0,k}\in L(E\cap(a^i)^\bot))
\]
of the summands $q^k$ themselves. We obtain
\[
\sum_{q^k_r>0}q^k_r\cdot t\left(p(r)\right) + t(q^{0,k}) = t\left(
\sum_{q^k_r>0}q^k_r\,p(r)+q^{0,k}\right) =: t(q^{1,k})
\]
and
\[
\sum_{q^k_r<0}q^k_r\cdot
t\left(p(r)\right)= t\left( \sum_{q^k_r<0}q^k_r\,p(r)\right)=:t(q^{2,k})\,.
\]
Up to elements of $E\cap (a^i)^\bot$, the relations $q^{1,k}$ and $q^{2,k}$ are
connected by
the common
\[
(q^{1,k})^-=-q^{1,k}_{r(i)}\cdot e^{r(i)}=\langle a^i,\bar{q}^k\rangle
\cdot e^{r(i)}=q^{2,k}_{r(i)}\cdot e^{r(i)}=(q^{2,k})^+\,.
\]
Hence, Lemma \zitat{5}{3} yields
\[
\sum_r q^k_r\cdot t\left(p(r)\right) + t(q^{0,k}) = t(q^{1,k}) + t(q^{2,k}) =
t\left(
q^{1,k}+q^{2,k}\right) = t(q^k)\,,
\]
and we conclude
\bean
\sum_k t(q^k) &=&
\sum_k \left(\sum_r q^k_r\cdot t\left(p(r)\right) + t(q^{0,k})\right)\\
&=&
\sum_r \left( \sum_k q^k_r \right) t\left(p(r)\right) + t\left(\sum_k
q^{0,k}\right)
\quad (\mbox{cf. previous claim})\\
&=&
0+ t\left( \sum_k q^k - \sum_{k,r} q^k_r\,p(r) \right)\\
&=& 0\,.
\vspace{2ex}
\eean
\par

Finally, only (c) is left. Let $a^i$, $a^j$ be two adjacent edges of $\sigma$.
We adapt the construction of the elementary relations
$p(r)$. Instead of the $r(i)$'s, we will use elements $r(i,j)\in E$
characterized by the
property
\[
\langle a^i, r(i,j)\rangle = 1\,,\; \langle a^j, r(i,j)\rangle = 0\,.
\]
(Those elements exist, since $Y$ is assumed to be smooth in codimension 2.)\\
Now, we define
\[
p(r):= e^r - \langle a^i,r\rangle \cdot e^{r(i,j)} - \langle a^j,r \rangle
\cdot e^{r(j,i)}
+ [\mbox{suitable element of }\Z^{E\cap(a^i)^\bot\cap(a^j)^\bot}]\,.
\]
These special $p(r)$'s meet the usual properties (i) and (ii) - but for the two
different
indices $i$ and $j$ at the same time. In particular, if $q\in L(E)$ is any
relation, then
$q$ and $\sum_{r\in E}q_r\cdot p(r)$ differ by some element of
$L(E\cap(a^i)^\bot\cap(a^j)^\bot)$ only.\\
\par

%
%
\sect{An alternative to the complex $L(E^R)_{\bullet}$}\label{s6}


\neu{61}
Let $R\in M$ be fixed for the whole \S \ref{s6}. The complex $L(E^R)_{\bullet}$
introduced in \zitat{3}{2} fits naturally into the exact sequence
\[
0\rightarrow L(E^R)_{\bullet} \longrightarrow (\Z^{E^R})_{\bullet}
\longrightarrow \mbox{span}(E^R)_{\bullet}\rightarrow 0
\]
of complexes built in the same way as $L(E^R)_{\bullet}$, i.e.
\[
(\Z^{E^R})_{-k} := \oplus\!\!\!\!\!\!_{\begin{array}{c}
\ks\tau<\sigma\vspace{-1ex} \\ \ks dim\, \tau=k \end{array}}
\!\!\!\!\Z^{E^R_{\tau}}
\qquad \mbox{and}\qquad
\mbox{span}(E^R)_{-k} := \oplus\!\!\!\!\!\!_{\begin{array}{c}
\ks\tau<\sigma\vspace{-1ex} \\ \ks dim\, \tau=k \end{array}}
\!\!\!\!\mbox{span}(E^R_{\tau})\,.
\]
\par

{\bf Lemma:}{\em
The complex $(\Z^{E^R})_{\bullet}$ is exact.\\
}
\par

{\bf Proof:}
The complex $(\Z^{E^R})_{\bullet}$ can be decomposed into a direct sum
\[
(\Z^{E^R})_{\bullet} = \bigoplus_{r\in M} (\Z^{E^R})(r)_{\bullet}
\]
showing the contribution of each $r\in M$. The complexes occuring as summands
are
defined as
\bean
(\Z^{E^R})(r)_{-k} &:=&
\oplus\!\!\!\!\!\!_{\begin{array}{c}
\ks\tau<\sigma\vspace{-1ex}\\ \ks dim\, \tau=k \end{array}} \!\!\!\!
\left\{ \begin{array}{ll}
\Z=\Z^{\{r\}} & \mbox{ for } r\in E^R_{\tau}\\
0 & \mbox{ otherwise}
\end{array} \right\}\\
&=&
\Z^{\#\{\tau\,|\; dim\,\tau=k; \, r\in E^R_{\tau}\}}\,.
\eean
Denote by $H^+$ the halfspace
\[
H^+ := \{ a\in N_{\R}\,|\; \langle a,r\rangle < \langle a, R\rangle\} \subseteq
N_{\R}.
\]
Then, for $\tau \neq 0$, the fact that $r\in E^R_{\tau}$ is equivalent to
$\tau \setminus \{0\} \subseteq H^+$. On the other hand, $r\in E^R_0$
corresponds to
the condition $\sigma \cap H^+ \neq \emptyset$.\\
In particular, $(\Z^{E^R})(r)_{\bullet}$,
shifted by one place, equals the complex for computing the reduced homolgy of
the
topological space $\cup \{\tau\,|\;\tau \setminus \{0\} \subseteq H^+\}
\subseteq \sigma$ cut
by some affine hyperplane. Since this space is contractable, the complex is
exact.
\hfill$\Box$\\
\par

{\bf Corollary:}{\em
The complexes $L(E^R)_{\bullet}^\ast$ and $\mbox{span}(E^R)_{\bullet}^\ast[1]$
are
quasiisomorphic. In particular, under the usual assumptions (cf. Theorem
\zitat{3}{3}), we obtain
\[
T^i_Y(-R) = H^i\left( \mbox{span}(E^R)_{\bullet}^\ast\otimes _{\Z}\C\right)\,.
\vspace{2ex}
\]
}
\par


\neu{62}
We define the $\R$-vector spaces
\bean
V^R_i &:= &\mbox{span}_{\R}(E_i^R)=\left\{
\begin{array}{l@{\quad\mbox{for}\;\:}l}
0 &\langle a^i,R\rangle\le 0\\
\left[ a^i=0\right] \subseteq M_{\R} & \langle a^i,R\rangle =1\\
M_{\R}=\R^n & \langle a^i,R\rangle\ge2
\end{array}
\right.\\
&& \qquad\qquad\qquad\qquad\qquad\qquad\qquad\qquad\qquad (i=1,\ldots,N),\;
\mbox{ and}\\
V^R_{\tau} &:=& \cap_{a^i\in \tau} V^R_i
\supseteq \mbox{span}_{\R}(E^R_{\tau})
\quad (\mbox{for faces } \tau<\sigma)\,.
\eean
$\,$
\vspace{-2ex}\\
\par

{\bf Proposition:}{\em
With ${\cal V}^R_{-k}:= \oplus\!\!\!\!\!\!_{\begin{array}{c}
\ks\tau<\sigma\vspace{-1.5ex} \\ \ks dim\, \tau=k \end{array}}
\!\!\!\!V_{\tau}^R$ we obtain a complex
${\cal V}^R_{\bullet}
\supseteq \mbox{span}_{\R}(E^R)_{\bullet}$.
Moreover, if $Y$ is smooth in
codimension $k$, then both complexes are equal at $\geq\!(-k)$.
}
\\
\par

{\bf Proof:}
$V^R_{\tau} = \mbox{span}_{\R}(E^R_{\tau})$ is true
for smooth cones $\tau<\sigma$ (cf.(3.7) of \cite{T1}).
\hfill$\Box$\\
\par

{\bf Corollary:}{\em
\begin{itemize}
\item[(1)]
If $Y$ is smooth in codimension 2, then $T^1_Y(-R)=
H^1\left(({\cal V}^R_{\bullet})^\ast \otimes_{\R} \C \right)$.
\item[(2)]
If $Y$ is smooth in codimension 3, then $T^2_Y(-R)=
H^2\left(({\cal V}^R_{\bullet})^\ast \otimes_{\R} \C \right)$.
\vspace{1ex}
\end{itemize}
}
\par

The formula (1) for $T^1_Y$ (with a more boring proof) was already obtained
in (4.4) of \cite{T1}.\\
\par

%
%
\sect{3-dimensional Gorenstein singularities}\label{s7}


\neu{71}
We want to apply the previous results for the special case of an isolated,
3-dimensional,
toric Gorenstein singularity. We start with fixing the notations.\\
\par
Let $Q=\mbox{conv}(a^1,\dots,a^N)\subseteq \R^2$ be a lattice polygon with
primitive
edges
\[
d^i:= a^{i+1}-a^i\in \Z^2\,.
\]
Embedding $\R^2$ as the affine hyperplane $[a_3=1]$ into $N_{\R}:=\R^3$,
we can define the cone
\[
\sigma:= \mbox{Cone}(Q) \subseteq N_{\R}\,.
\]
The fundamental generators of $\sigma$ equal the vectors
$(a^1,1),\dots,(a^N,1)$, which we
will also denote by $a^1,\dots,a^N$, respectively. \\
\par
The vector space $M_{\R}$ contains a special element $R^\ast:=[0,0;1]$:
\begin{itemize}
\item
$\langle \bullet,R^\ast\rangle = 1$ defines the affine hyperplane containing
$Q$,
\item
$\langle \bullet,R^\ast\rangle = 0$ describes the vectorspace containing the
edges
$d^i$ of $Q$.
\end{itemize}
The structure of the dual cone $\sigma^{\veee}$ can be described as follows:
\begin{itemize}
\item
$[c;\eta]\in M_{\R}$ is contained in $\sigma^{\veee}$, iff $\langle Q,-c\rangle
\leq \eta$.
\item
$[c;\eta]\in \partial \sigma^{\veee}$ iff there exists some $i$ with
$\langle a^i,-c\rangle = \eta$.
\item
The set $E$ contains $R^\ast$. However, $E\setminus \{R^\ast\}\subseteq
\partial
\sigma^{\veee}$.
\vspace{1ex}
\end{itemize}
\par

{\bf Remark:}
The toric variety $Y$ built by the cone $\sigma$ is 3-dimensional, Gorenstein,
and
regular outside its 0-dimensional orbit. Moreover, all those singularities can
be
obtained in this way.\\
\par


\neu{72}
Let $V$ denote the $(N-2)$-dimensional $\R$-vector space
\[
V:=\{(t_1,\dots,t_N)\,|\; \sum_i t_i\,d^i=0\}\subseteq \R^N\,.
\]
The non-negative tuples among the $\ku{t}\in V$ describe the set of Minkowski
summands
$Q_{\ku{t}}$
of positive multiples of the polygon $Q$. ($t_i$ is the scalar by which $d^i$
has to be
multiplied to get the $i$-th edge of $Q_{\ku{t}}$.)\\
\par
We consider the bilinear map
\[
\begin{array}{cclcl}
V&\times& \R^E & \stackrel{\Psi}{\longrightarrow}& \R\\
\ku{t}&,&[c;\eta]\in E &\mapsto&
\left\{ \begin{array}{ll}
0& \mbox{ if } c=0\quad(\mbox{i.e. } [c;\eta]=R^\ast)\\
\sum_{v=1}^{i-1} t_v\cdot
\langle d^v,-c\rangle & \mbox{ if }\langle a^i,-c\rangle =\eta\,.
\end{array} \right.
\end{array}
\]
Assuming both $a^1$ and the associated vertices of all Minkowski sumands
$Q_{\ku{t}}$
to coincide with $0\in\R^2$, the map $\Psi$  detects the maximal values of the
linear
functions $c$ on these summands
\[
\Psi(\ku{t},[c;\eta]) = \mbox{Max}\,(\langle a,-c\rangle\,|\; a\in
Q_{\ku{t}})\,.
\]
In particular, $\Psi(\ku{1},[c;\eta])=\eta$, i.e. $\Psi$ induces a map
\[
\Psi: \quad^{\kd V}\!\!/\!_{\kd \R\cdot\ku{1}} \times L_{\R}(E) \longrightarrow
\R\,.
\]
The results of \cite{Gor} and \cite{Sm} imply that $\Psi$ provides an
isomorphism
\[
^{\kd V_{\C}}\!\!/\!_{\kd \C\cdot\ku{1}}\stackrel{\sim}{\longrightarrow}
\left(\left. ^{\displaystyle L(E_0^{R^\ast})}\!\!\right/
\!_{\displaystyle \sum_i L(E_i^{R^\ast})} \right)^\ast \otimes_{\Z} \C
\cong T^1_Y(-R^\ast)= T^1_Y\,.
\]
In particular, $\mbox{dim}\,T^1_Y = N-3$.\\
\par


\neu{73}
Let $R\in M$. Combining the general results of \S \ref{s6} with the fact
\[
\bigcap_i V_i^R = \mbox{Ker}\left[ \oplus_i(V_i^R\cap V^R_{i+1})\longrightarrow
\oplus_i V^R_i\right]
\]
coming from the special situation we are in, we obtain the handsome formula
\[
T^2_Y(-R)=
\left[ \left.^{\kd \bigcap_i (\mbox{span}_{\C} E_i^R)}\!\!\! \right/
\!\!\! _{\kd \mbox{span}_{\C} (\bigcap_i E_i^R)} \right] ^\ast\,.
\]
$T^1_Y$ is concentrated in the degree $-R^\ast$. Hence, for computing $T^2_Y$,
the
degrees $-kR^\ast$ ($k\geq 2$) are the most interesting (but not only) ones. In
this
special case, the vector spaces $V^{kR^\ast}_i$ equal $M_{\R}$, i.e.
\[
T^2_Y(-kR^\ast)= \left[ \left. ^{\kd M_{\C}}\!\!\! \right/ \!\!\!
_{\kd \mbox{span}_{\C} (\bigcap_i E_i^{kR^\ast})} \right] ^\ast
\subseteq
\left[ \left. ^{\kd M_{\C}}\!\!\! \right/ \!\!\! _{\kd \C\cdot R^\ast}\right]
^\ast =
\mbox{span}_{\C}(d^1,\dots,d^N) \subseteq N_{\C}\,.
\vspace{2ex}
\]
\par

{\bf Proposition:}
{\em
For $c\in \R^2$ denote by
\[
d(c):= \mbox{Max}\,(\langle a^i,c\rangle\,|\; i=1,\dots,N) -
\mbox{Min}\,(\langle a^i,c\rangle\,|\; i=1,\dots,N)
\]
the diameter of $Q$ in $c$-direction. If
\[
k_1:= \!\!\begin{array}[t]{c}
\mbox{Min}
\vspace{-1ex}\\
\ks c\in\Z^2\setminus 0
\end{array} \!\!d(c) \quad \mbox{ and } \quad
k_2:= \!\!\!\begin{array}[t]{c}
\mbox{Min}
\vspace{-1ex}\\ \ks c,c'\in\Z^2
\vspace{-1ex}\\ \ks lin.\, indept.
\end{array} \!\!\!\mbox{Max}\,[ d(c), d(c')]\,,
\]
then
$\quad\begin{array}[t]{lll}
\dim T^2_Y(-kR^\ast) = 2 & \mbox{ for } & 2\leq k \leq k_1\,,\\
\dim T^2_Y(-kR^\ast) = 1 & \mbox{ for } & k_1+1\leq k \leq k_2\,,\mbox{ and}\\
\dim T^2_Y(-kR^\ast) = 0 & \mbox{ for } & k_2+1\leq k \,.
\end{array}
\vspace{2ex}
$
}
\par

{\bf Proof:}
We have to determine the dimension of $\;\mbox{span}_{\C}\left( \bigcap_i
E_i^{kR^\ast}\right)\!\!\left/\!\!_{\kd \C\cdot R^\ast}\right.$. Computing
modulo $R^\ast$
simply means to forget the $\eta$ in $[c;\eta]\in M$. Hence, we are done by the
following
observation for each $c\in \Z^2\setminus0$:
\[
\begin{array}{rcl}
\exists \eta\in\Z: \;[c,\eta]\in \bigcap_iK_i^{kR^\ast}
& \Longleftrightarrow &
\exists \eta\in\Z: \;(k-1)R^{\ast}\kgeq [c;\eta] \kgeq 0\\
& \Longleftrightarrow &
d(c) \leq k-1\,.
\end{array}
\vspace{-3ex}
\]
\hfill$\Box$\\
\par

{\bf Corollary:}
{\em
Unless $Y=\C^3$ or $Y=\mbox{cone over }\PP^1\times\PP^1$, we have
\[
T^2_Y(-2R^\ast)= \mbox{span}_{\C}(d^1,\dots,d^N),
\]
i.e. $\dim T^2_Y(-2R^\ast)=2$.}
\\
\par


\neu{74}
{\bf Proposition:}
{\em
Using both the isomorphism $\;V_{\C}\!\!\left/\!\!_{\kd \C\cdot\ku{1}}\right.
\stackrel{\sim}{\rightarrow} T^1_Y$ and the injection
$T^2_Y(-2R^\ast)\hookrightarrow \mbox{span}_{\C}(d^1,\dots,d^N)$,
the cup product $T^1_Y\times T^1_Y \rightarrow T^2_Y$ equals the bilinear map
\[
\begin{array}{ccccc}
V_{\C}\!\!\left/\!\!_{\kd \C\cdot\ku{1}}\right. &
\times &
V_{\C}\!\!\left/\!\!_{\kd \C\cdot\ku{1}}\right. &
\longrightarrow &
\mbox{span}_{\C}(d^1,\dots,d^N)\\
\ku{s} &
,
&
\ku{t} &
\mapsto &
\sum_i s_i\,t_i\,d^i\,.
\end{array}
\vspace{1ex}
\]
}
\par

{\bf Proof:}
{\em Step 1:} To apply Theorem \zitat{3}{6} we will combine the isomorphisms
for $T^2_Y$ presented in \S \ref{s6} and \zitat{7}{3}. Actually, we will
describe the dual map
by associating to each $r\in M$ an element
$[q^1(r),\dots,q^N(r)]\in\oplus_iL(E_i^{2R^\ast})$.\\
First, for every $i=1,\dots,N$, we have to write
$r\in M = (\mbox{span}\, E_i^{2R^\ast})\cap (\mbox{span}\, E_{i+1}^{2R^\ast})$
as a linear
combination of elements from $E_i^{2R^\ast}\cap E_{i+1}^{2R^\ast}$.
This set contains a $\Z$-basis for $M$ consisting of
\begin{itemize}
\item
$r^i:=$ primitive element of $\sigma^{\veee}\cap (a^i)^\bot \cap
(a^{i+1})^\bot$,
\item
$R^\ast$, and
\item
$r(i):= r(i,i+1)$ (cf. notation at the end of \zitat{5}{4}), i.e.
$\begin{array}[t]{l}
\langle a^i, r(i)\rangle = 1 \mbox{ and}\\
\langle a^{i+1}, r(i) \rangle = 0\,.
\end{array} $
\end{itemize}
In particular, we can write
\[
r= g^i(r)\cdot r^i + \langle a^{i+1},r\rangle\cdot R^\ast + \left(
\langle a^i,r \rangle - \langle a^{i+1},r\rangle \right) \cdot r(i)
\]
with some integer $g^i(r)\in \Z$.\\
\par
Now, we have to apply the differential in the complex
$(\Z^{E^{2R^\ast}})_{\bullet}$,
i.e. we map the previous expression via the map
\[
\oplus_i \Z^{E_i^{2R^\ast}\cap E_{i+1}^{2R^\ast}} \longrightarrow
\oplus_i \Z^{E_i^{2R^\ast}}\, .
\]
The result is (for every $i$) the element of $L(E_i^{2R^\ast})$
\[
\begin{array}{l}
g^i(r)\, e^{r^i} - g^{i-1}(r)\, e^{r^{i-1}} + \langle a^i-a^{i+1},r\rangle\cdot
e^{r(i)} -
\langle a^{i-1} - a^i,r\rangle \cdot e^{r(i-1)} +
\langle a^{i+1}-a^i, r \rangle\cdot e^{R^\ast}
\vspace{1ex}\\
\qquad = \langle d^i,r\rangle \cdot \left( e^{R^\ast} -
e^{r(i)}\right) + [(a^i)^\bot \mbox{-summands}] =: q^i(r)\,.
\end{array}
\vspace{2ex}
\]
\par

{\em Step 2:}
Defining
\[
q^i:= e^{R^\ast}-e^{r(i)} + [(a^i)^\bot \mbox{-summands}] \in
L(E_i^{2R^\ast})\quad
(i=1,\dots,N)\,,
\]
we use Theorem \zitat{3}{6} and the second remark of \zitat{3}{6} to obtain
\[
(\ku{s}\cup\ku{t})_i \left( q^i(r) \right) =
\langle d^i,r \rangle \cdot t_{\Psi(\ku{s},\bullet), \Psi(\ku{t},\bullet),
R^\ast,R^\ast}
(q^i) = \Psi(\ku{s},q^i)\cdot \Psi(\ku{t},q^i)\,.
\]
To compute those two factors, we take a closer look at the $q^i$'s. Let
\[
q^i= e^{R^\ast}-e^{r(i)} + \sum_v \lambda^i_v \,e^{[c^v; \eta^v]}\,,
\vspace{-1ex}
\]
and the sum is taken over those $v$'s meeting the property
$\langle a^i,-c^v\rangle = \eta^v$. Then, by definition of $\Psi$ in
\zitat{7}{2},
\[
\Psi(\ku{s},q^i)= \sum_{j=1}^{(i+1)-1}s_j\,\langle d^j, r(i)\rangle -
\sum_v \lambda^i_v \cdot \left( \sum_{j=1}^{i-1}s_j\, \langle d^j,c^v\rangle
\right)\,.
\]
On the other hand, we know that $q^i$ is a relation, i.e. the equation
\[
R^\ast - r(i) + \sum_v \lambda^i_v [c^v; \eta^v] =0
\vspace{-1ex}
\]
is true in $M$. Hence,
\[
\begin{array}{rcl}
\Psi(\ku{s},q^i)&=& \sum_{j=1}^i s_j\,\langle d^j, r(i)\rangle -
\sum_{j=1}^{i-1} s_j \langle d^j, r(i)\rangle
\vspace{0.5ex}\\
&=& s_i\cdot \langle d^i, r(i)\rangle
\vspace{0.5ex}\\
&=& -s_i\,.
\end{array}
\vspace{-3ex}
\]
\hfill$\Box$
\vspace{2ex}\\
\par

$T^1_Y\subseteq \C^N$ is the tangent space of the versal base space $S$ of our
singularity
$Y$. It is given by the linear equation $\sum_i t_i\cdot d^i=0$.\\
On the other hand,
the cup product $T^1_Y\times T^1_Y\rightarrow T^2_Y$ shows the quadratic part
of the
equations defining $S\subseteq\C^N$. By the previous proposition, it equals
$\sum_i t_i^2 \cdot d^i$.\\
\par
These facts suggest an idea how the equations of $S\subseteq\C^N$ could look
like. In
\cite{Vers} we have proved this conjecture; $S$ is indeed given by the
equations
\[
\sum_{i=1}^N t_i^k \cdot d^i =0\quad (k\geq 1)\,.
\vspace{3ex}
\]
\par

%
%

\end{document}